\documentclass[12pt]{article}
\usepackage{amsmath,amsfonts,amssymb}

\newcommand{\Span}[1]{\left<#1\right>}
\newcommand{\fie}{\varphi}
\newcommand{\al}{\alpha}
\newcommand{\be}{\beta}
\newcommand{\ga}{\gamma}
\newcommand{\de}{\delta}

\newcommand{\la}{\lambda}
\newcommand{\si}{\sigma}

\newcommand{\Ga}{\Gamma}

\newcommand{\Si}{\Sigma}

\newcommand{\C}{\mathbb C}
\newcommand{\FF}{\mathbb F}
\newcommand{\PP}{\mathbb P}
\newcommand{\wP}{w\mathbb P}
\newcommand{\Q}{\mathbb Q}
\newcommand{\N}{\mathbb N}
\newcommand{\Z}{\mathbb Z}
\newcommand{\Oh}{\mathcal O}
\newcommand{\sB}{\mathcal B}

\newcommand{\sI}{\mathcal I}
\newcommand{\sJ}{\mathcal J}

\newcommand{\sN}{\mathcal N}
\newcommand{\bara}{\overline a}
\newcommand{\ba}{\mathbf a}
\newcommand{\bb}{\mathbf b}
\newcommand{\bx}{\mathbf x}
\newcommand{\by}{\mathbf y}
\newcommand{\broken}{\dasharrow}
\newcommand{\ot}{\leftarrow}

\newcommand{\onto}{\twoheadrightarrow}
\DeclareMathOperator{\rank}{rank}
\DeclareMathOperator{\wt}{wt}
\DeclareMathOperator{\Cl}{Cl}
\DeclareMathOperator{\Grass}{Grass}

\DeclareMathOperator{\Jac}{Jac}
\DeclareMathOperator{\Pf}{Pf}

\DeclareMathOperator{\Rad}{Rad}
\DeclareMathOperator{\Sing}{Sing}

\newtheorem{thm}{Theorem}[section]
\newtheorem{lem}[thm]{Lemma}
\newtheorem{exa}[thm]{Example}
\newtheorem{exun}{Example}

\newtheorem{dfn}[thm]{Definition}

\newcommand{\QED}{\ifhmode\unskip\nobreak\fi\quad\ensuremath{\mathrm{QED}}}

\numberwithin{equation}{section}

\begin{document}

\title{Fano 3-folds in codimension 4, \\ Tom and Jerry. Part I}
\author{Gavin Brown \and Michael Kerber \and Miles Reid}
\date{}

\maketitle

\begin{abstract}
This work is part of the Graded Ring Database project \cite{GRDB}, and
is a sequel to \cite{APhD} and \cite{ABR}. We introduce a strategy based
on Kustin--Miller unprojection that allows us to construct many hundreds
of Gorenstein codimension~4 ideals with \hbox{$9\times 16$} resolutions (that
is, 9 equations and 16 first syzygies). Our two basic games are called
Tom and Jerry; the main application is the biregular construction of
most of the anticanonically polarised Mori Fano 3-folds of Alt{\i}nok's
thesis \cite{APhD}. There are 115 cases whose numerical data (in effect,
the Hilbert series) allow a Type~I projection. In every case, at least
one Tom and one Jerry construction works, providing at least two
deformation families of quasismooth Fano \hbox{3-folds} having the same
numerics but different topology.



\end{abstract}

\setcounter{tocdepth}{1}

\section{Introduction and the classification of\\Fano 3-folds}\label{s!intr}

A {\em Fano $3$-fold} $X$ is a normal projective 3-fold whose
anti\-canonical divisor $-K_X=A$ is \hbox{$\Q$-Cartier} and ample. We
eventually impose additional conditions on the singularities and class
group of~$X$, such as terminal, $\Q$-factorial, quasi\-smooth, prime
(that is, class group $\Cl X$ of rank~1) or $\Cl X=\Z\cdot A$, but more
general cases occur in the course of our arguments.

We study $X$ via its anticanonical graded ring
\[
R(X,A) = \bigoplus_{m \in \N} H^0(X,mA).
\]
Choosing generators of $R(X,A)$ embeds $X$ as a projectively normal
sub\-variety $X\subset\PP(a_1,\dots,a_n)$ in weighted projective space.
The anti\-canonical ring $R(X,A)$ is known to be Gorenstein, and we say
that $X\subset\PP(a_1,\dots,a_n)$ is projectively Gorenstein. The {\em
codimension} of $X$ refers to this anticanonical embedding. The
discrete invariants of a Fano 3-fold $X$ are its {\em genus} $g$
(defined by $g+2=h^0(X,-K_X)$) together with a basket of terminal cyclic
quotient singularities; for details see \ref{s!num} and \cite{ABR}.

In small codimension we can write down hypersurfaces, codimension~2
complete intersections and codimension~3 Pfaffian varieties fluently.
This underlies the classification of Fano $3$-folds in
codimension~$\le3$ (see \cite{C3f}, \cite{Fl} and \cite{APhD}): the
famous 95 weighted hypersurfaces, 85 codimension~2 families, and 70
families in codimension~3, of which 69 are $5\times5$ Pfaffians.
Gorenstein in codimension~4 remains one of the frontiers of science:
there is no automatic structure theory, and deformations are almost
always obstructed. Type~I projection and Kustin--Miller unprojection
(see \cite{KM}, \cite{PR}, \cite{Ki}) is a substitute that is sometimes
adequate. This paper addresses codimension~4 Fano $3$-folds in this vein.

The analysis of \cite{APhD}, \cite{A}, \cite{ABR}, \cite{GRDB} provides
145 numerical candidates for codimension~4 Fano $3$-folds. This paper
isolates 115 of these that can be studied using Type~I projections,
hence as Kustin--Miller unprojections. Our main result is
Theorem~\ref{th!main}: each of these 115 numerical candidates occurs in
at least two ways (the Tom and Jerry of the title), that give rise to
topologically distinct varieties $X$.

The reducibility of the Hilbert scheme of Fano 3-folds is a systematic
feature of our results, that goes back to Takagi's study of prime Fano
3-folds with basket of $\frac12(1,1,1)$ points (\cite{T}, Theorem~0.3).
He describes families of varieties having the same invariants, but
arising from different ``Takeuchi programs'', that is, different
Sarkisov links. Four of his numerical cases have codimension~4. The
first, No.~1.4 in the tables of \cite{T}, is our initial case
$X\subset\PP^7(1^7,2)$; it projects to the $(2,2,2)$ complete
intersection, so has $7\times12$ resolution and is unrelated to Tom and
Jerry. Takagi's three other pairs of codimension~4 cases correspond to
our Tom and Jerry families as follows:
\begin{equation*}
\begin{array}{ccc}
X\subset\PP^7(1^4,2^4): & \hbox{Tom$_1 =$ No.~2.2 (8 nodes)}, & \hbox{Jer$_{45} =$ No.~3.3 (9 nodes)},\\
X\subset\PP^7(1^5,2^3): & \hbox{Tom$_1 =$ No.~5.4 (7 nodes)}, & \hbox{Jer$_{23} =$ No.~4.1 (8 nodes)},\\
X\subset\PP^7(1^6,2^2): & \hbox{Tom$_1 =$ No.~4.4 (6 nodes)}, & \hbox{Jer$_{15} =$ No.~1.1 (7 nodes)}.
\end{array}
\end{equation*}
Each of these is prime. In our treatment, each of these numerical cases
admits one further Jerry family consisting of Fano 3-folds of Picard rank
$\ge2$.

Section~\ref{sec!anc} traces the origin of Tom and Jerry back to the
geometry of linear subspaces of $\Grass(2,5)$ and associated
unprojections to twisted forms of $\PP^2\times\PP^2$ and
$\PP^1\times\PP^1\times\PP^1$; for more on this, see
Section~\ref{s!fmt}. Section~\ref{s!main} is a detailed discussion of
our Main Theorem~\ref{th!main}, whose proof occupies the rest of the
paper. Flowchart~\ref{ss!flchart} maps out the proof, which involves
many thousand computer algebra calculations. Section~\ref{s!fmt}
discusses the wider issue of codimension~4 formats, and serves as a
mathematical counterpart to the computer algebra of
Sections~\ref{s!fail}--\ref{s!grdb}. We do not elaborate on this point,
but Tom and Jerry star in many other parallel or serial unprojection
stories beyond Fano 3-folds or codimension~4, notably the diptych
varieties of \cite{aflip}.

We are indebted to a referee for several pertinent remarks that led to
improvements, and to a second referee who verified our computer algebra
calculations independently. This research is supported by the Korean
government WCU Grant R33-2008-000-10101-0.

\section{Ancestral examples}\label{sec!anc}

\subsection{Linear subspaces of $\Grass(2,5)$}\label{s!dP6}

A del Pezzo variety of degree~5 is an $n$-fold $Y^n_5\subset\PP^{n+3}$
of codimension~3, defined by 5 quadrics that are Pfaffians of a
$5\times5$ skew matrix of linear forms. Thus $Y$ is a linear section of
Pl\"ucker $\Grass(2,5)\subset\PP(\bigwedge^2V)$ (here $V=\C^5$). We want
to unproject a projective linear subspace $\PP^{n-1}$ contained as a
divisor in $Y$ to construct a degree~6 del Pezzo variety
$X^n_6\subset\PP^{n+4}$. The crucial point is the following.

\begin{lem} The Pl\"ucker embedding $\Grass(2,5)$ contains two families
of maximal linear subspaces. These arise from
\begin{enumerate}
\renewcommand{\labelenumi}{(\Roman{enumi})}

\item The $4$-dimensional vector subspace $v\wedge V\subset\bigwedge^2V$
for a fixed $v\in V$.

\item The $3$-dimensional subspace $\bigwedge^2U\subset\bigwedge^2V$ for
a fixed $3$-dimensional vector subspace $U\subset V$.
\end{enumerate}
\end{lem}

Thus there are two different formats to set up $\PP^{n-1}\subset Y$.
Case~I gives $\PP^3_v\subset\Grass(2,5)$. A section of $\Grass(2,5)$ by
a general $\PP^7$ containing $\PP^3_v$ is a 4-fold $Y^4$ whose
unprojection is $\PP^2\times\PP^2\subset\PP^8$. Case~II gives
$\Grass(2,U)=\PP^2_U\subset\Grass(2,5)$. A section of $\Grass(2,5)$ by a
general $\PP^6$ containing $\PP^2_U$ is a 3-fold $Y^3$ whose unprojection
is $\PP^1\times\PP^1\times\PP^1\subset\PP^7$.

The proof is a lovely exercise. Hint: use local and Pl\"ucker coordinates
\begin{equation}
\begin{pmatrix} 1&0&a_1&a_2&a_3 \\ 0&1&b_1&b_2&b_3 \end{pmatrix}
\quad\hbox{and}\quad
\begin{pmatrix} 1&a_1&a_2&a_3 \\ &b_1&b_2&b_3 \\ &&m_{12} & m_{13} \\
&&& m_{23} \end{pmatrix}
\end{equation}
with Pl\"ucker equations $m_{12}=a_1b_2-a_2b_1$, etc.; permute the
indices and choose signs pragmatically to make this true. Prove that in
Pl\"ucker $\PP^9$, the tangent plane $m_{12}=m_{13}=m_{23}=0$ intersects
$\Grass(2,5)$ in the cone over the Segre embedding of $\PP^1\times\PP^2$.

\subsection{Tom$_1$ and Jer$_{12}$ in equations}\label{s!same}
Tom$_1$ is
\begin{equation}
\begin{pmatrix}
y_1&y_2&y_3&y_4 \\
& m_{23}& m_{24}& m_{25}\\
&& m_{34}& m_{35}\\
&&& m_{45}\\
\end{pmatrix} \label{eq!1.2}
\end{equation}
with $y_{1\dots4}$ arbitrary elements, and the six entries $m_{ij}$
linear combinations of a regular sequence $x_{1\dots4}$ of length four.
Expressed vaguely, there are ``two constraints on these six entries'';
these two coincidences take the simplest form when $m_{23}=m_{45}=0$. In
this case, the Pfaffian equations all reduce to binomials, and can be
seen as the $2\times2$ minors of an array: as a slogan,
\begin{equation}
\hbox{$4\times4$ Pfaffians of }
\left(\begin{smallmatrix}
y_1&y_2&y_3&y_4 \\
& 0 & m_{24}& m_{25}\\
&& m_{34}& m_{35}\\
&&& 0\\
\end{smallmatrix}\right)
=\hbox{$2\times2$ minors of }
\left(\begin{smallmatrix}
* & y_3 &y_4 \\
y_1 & m_{24}& m_{25}\\
y_2& m_{34}& m_{35}
\end{smallmatrix}\right).
 \label{eq!1.4}
\end{equation}
That is, the $4\times 4$ Pfaffians on the left equal the five $2\times2$
minors of the array on the right. To see the Segre embedding of
$\PP^2\times\PP^2$ and its linear projection from a point, replace the
star entry by the unprojection variable $s$.

In a similar style, Jer$_{12}$ is
\begin{equation}
\begin{pmatrix}
m_{12} & m_{13}& m_{14}& m_{15}\\
& m_{23}& m_{24}& m_{25}\\
&& y_{34}& y_{35}\\
&&& y_{45}\\
\end{pmatrix}
\end{equation}
with $y_{34},y_{35},y_{45}$ arbitrary, and the seven entries $m_{ij}$
linear combinations of $x_{1\dots4}$. Vaguely, ``three constraints on
these seven entries''; most simply, these take the form
$m_{15}=m_{23}=0$, $m_{24}=m_{14}$. We leave you to see this as the
linear projection of $\PP^1\times\PP^1\times\PP^1$, starting from the
hint:
\begin{equation} \label{eq!cube1}
\hbox{$4\times4$ Pfaffians of }
\begin{pmatrix}
t&z_1&z_2&0 \\
& 0 & z_2& z_3\\
&& y_3& y_2\\
&&& y_1\\
\end{pmatrix}
=\hbox{$2\times2$ minors of }
\begin{picture}(30,30)(0,0)
\renewcommand{\arraycolsep}{.2em}
\put(8,-7){$\begin{matrix}
* &\frac{\quad}{\quad}& y_2 \\
\vert&&\vert \\
y_1 &\frac{\quad}{\quad}& z_3
\end{matrix}$}
\put(28,7){$\begin{matrix}
y_3 &\frac{\quad}{\quad}&z_1 \\
\vert&&\vert \\
z_2 &\frac{\quad}{\quad}& t
\end{matrix}$}
\qbezier(20,15)(23,17)(26,19)
\qbezier(20,-12)(23,-10)(26,-8)
\qbezier(52,15)(55,17)(58,19)
\qbezier(52,-12)(55,-10)(58,-8)
\end{picture}
\kern16mm
\end{equation}
that is, on the right, take $2\times 2$ minors of the three square faces
out of $t$, together with the ``diagonal'' minors $y_1z_1=y_2z_2=y_3z_3$,
then replace the star by an unprojection variable. Compare
\eqref{eq!cube}.

\subsection{General conclusions}
\begin{dfn} \rm \label{d!TJ}
Tom$_i$ and Jer$_{ij}$ are matrix formats that specify unprojection
data, namely a codimension~3 scheme $Y$ defined by a $5\times5$ Pfaffian
ideal, containing a codimension~4 complete intersection $D$. Given a
regular sequence $x_{1\dots4}$ in a regular ambient ring $R$ generating
the ideal $I_D$, the ideal of $Y$ is generated by the Pfaffians of a
$5\times5$ skew matrix $M$ with entries in $R$, subject to the
conditions \begin{description}

\item{Tom$_i$:} the 6 entries $m_{jk}\in I_D$ for all $j,k\ne i$; in
other words, the 4 entries $m_{ij}$ of the $i$th row and column are free
choices, but the other entries of $M$ are required to be in $I_D$. See
\eqref{eq!2.8} for an example.

\item{Jer$_{ij}$:} the 7 entries $m_{kl}\in I_D$ if either $k$ or $l$
equals $i$ or $j$. See \eqref{eq!2.9} for an example. The bound entries
are the {\em pivot} $m_{ij}$ and the two rows and columns through it.
The 3 free entries are the Pfaffian partners $m_{kl}$, $m_{km}$,
$m_{lm}$ of the pivot, where $\{i,j,k,l,m\}=\{1,2,3,4,5\}$. In $Y$, the
pivot vanishes twice on $D$.

\end{description}
\end{dfn}

Case~I in \ref{s!dP6} is the ancestor of our Tom constructions and~II
that of Jerry. Our main aim in what follows is to work out several
hundred applications of the same formalism to biregular models of Fano
3-folds, when our ``constraints''
\begin{equation}
m_{ij}=\hbox{linear combination of $x_{1\dots4}$}
\end{equation}
are not linear, do not necessarily reduce to a simple normal form, and
display a rich variety of colourful and occasionally complicated
behaviour.

Nevertheless, the same general tendencies recur again and again. Tom
tends to be fatter than Jerry. Jerry tends to have a singular locus of
bigger degree than Tom, and the unprojected varieties $X$ have different
topologies, in fact different Euler numbers. For example, $Y^4$ in
Case~I has two lines of transversal nodes; the $Y^3$ in Case~II has
three nodes. If we only look at 3-folds in \ref{s!dP6} (cutting $Y^4$ by
a hyperplane), the unprojected varieties $X$ are then the familiar del
Pezzo 3-folds of index 2, namely the flag manifold of $\PP^2$ versus
$\PP^1\times\PP^1\times\PP^1$; see Section~\ref{s!Jac} (especially
Remark~\ref{rk!2v3}) for the number of nodes (2 and 3 in the two cases)
via enumerative geometry. Tom equations often relate to extensions of
$\PP^2\times\PP^2$ such as the ``extrasymmetric $6\times6$ format'';
Jerry equations often relate to extensions of
$\PP^1\times\PP^1\times\PP^1$ such as the ``rolling factors format'' (an
anticanonical divisor in a scroll) or the ``double Jerry format'';
Section~\ref{s!fmt} gives a brief discussion.

\section{The main result}\label{s!main}

\subsection{Numerical data  of Fano 3-folds}\label{s!num}
Let $X$ be a Fano $3$-fold. As explained in \cite{ABR}, the numerical
data of $X$ consists of an integer genus $g\ge-2$ plus a basket
$\sB=\{\frac1r(1,a,r-a)\}$ of terminal cyclic orbi\-fold points; this
data determines the Hilbert series $P_X(t)= \sum_{a\ge0}h^0(X,nA) t^n$ of
$R(X,A)$, and is equivalent to it. At present we only treat cases when
the ring is generated as simply as possible, and not (say) cases that
fall in a monogonal or hyperelliptic special case. The database
\cite{GRDB} lists cases of small codimension, including 145 candidate
cases in codimension~4 from Alt{\i}nok's thesis \cite{APhD}. We
sometimes say Fano \hbox{3-fold} to mean numerical candidate; the abuse of
terminology is fairly harmless, because practically all the candidates
in codimension $\le5$ (possibly all of them) give rise to quasismooth
Fano \hbox{3-folds}; in fact usually more than one family, as we now
relate.

\subsection{Type~I centre and Type I projection}\label{s!TyIc}
An orbifold point $P\in X$ of type $\frac1r(1,a,r-a)$ with $r\ge2$ is a
{\em Type~I centre} if its orbinates are restrictions of global forms
$x\in H^0(A)$, $y\in H^0(aA)$, $z\in H^0((r-a)A)$ of the same weight.
The condition means that after projecting, the exceptional locus of the
projection is a weighted projective plane $\PP(1,a,r-a)$ that is
embedded projectively normally.

One may view a projection $P\in X\broken Y\supset D$ in simple terms: in
geometry, as the map
$(x_1,\dots,x_n)\mapsto(x_1,\dots,\widehat{x_i},\dots,x_n)$ analogous to
linear projection $\PP^n\broken\PP^{n-1}$ from centre
$P_i=(0,\dots,1,\dots,0)$; or in algebra, as eliminating a variable,
corresponding to passing to a graded subring
$k[x_1,\dots,\widehat{x_i},\dots,x_n]$; to be clear, the distinguishing
characteristic is not the eliminated variable $x_i$, rather the point
$P_i$ and the complementary system of variables $x_j$ that vanish there.

We take the more sophisticated view of \cite{CPR}, 2.6.3 of a projection
as an intrinsic biregular construction of Mori theory; namely a diagram
\begin{equation}
\renewcommand{\arraycolsep}{1pt}
\begin{array}{ll}
\kern2.5em P \in X & \subset \PP(a_0,\dots,a_n)\\
\kern3.25em\nearrow \\
E \subset X_1 \\
\kern3.25em\searrow \\
\kern2.5em D \subset Y & \subset \PP(a_0,\dots,\widehat{a_k},\dots,a_n)
\end{array}
\label{eq!proj}
\end{equation}
consisting of an extremal extraction $\si\colon X_1\to X$ centred at $P$
followed by the anticanonical morphism $\fie\colon X_1\to Y$. In more
detail, we have the following result.

\begin{lem}\label{l!s-amp}
The Type I assumption implies that $-K_{X_1}$ is semiample.
The anticanonical morphism
$\fie\colon X_1\to Y$ contracts only curves $C$ with $-K_{X_1}C=0$
meeting the exceptional divisor $E=\PP(1,a,r-a)\subset X_1$
transversely in one point.
\end{lem}

\paragraph{Proof} A theorem of Kawamata \cite{Ka} (discussed also in
\cite{CPR}, 3.4.2) says that the $(1,a,r-a)$ weighted blowup $\si\colon
X_1\to X$ is the unique Mori extremal extraction whose centre meets the
$\frac1r(1,a,r-a)$ orbifold point $P\in X$. It has exceptional divisor
the weighted plane $E=\PP(1,a,r-a)$ with discrepancy $\frac1r$. Thus $-K_{X_1}=-K_X-\frac1rE$, and the anticanonical ring of $X_1$ consists
of forms of weight $d$ in $R(X,-K_X)$ vanishing to order $\ge\frac dr$
on $E$. The homogenising variable $x_k$ of degree $r$
with $x_k(P)=1$ does not vanish at all, so is eliminated. By
assumption, the orbinates $x,y,z$ at $P$ are global forms of weights
$1,a,r-a$ vanishing to order exactly $\frac1r,\frac ar,\frac{r-a}r$, so
these extend to regular elements of $R(X_1,-K_{X_1})$. Locally at $P$,
appropriate monomials in $x,y,z$ base the sheaves $\Oh_X(d)$ modulo any
power of the maximal ideal $m_P$, so we can adjust the remaining
generators $x_l$ of $R(X,-K_X)$ to vanish to order $\ge \frac{\wt x_l}r$,
and so they lift to $R(X_1,-K_{X_1})$. It follows that $-K_{X_1}$ is
semi\-ample and the anti\-canonical morphism $\fie\colon X_1\to Y$
takes $E$ isomorphically to $D\subset Y$.
\QED \par\medskip

In our cases, $\fie$ contracts a nonempty finite set of flopping
curves to singular points of $Y$ on $D$, and $Y$ is a codimension~3
Fano 3-fold. The anti\-canonical model $Y$ is not $\Q$-factorial because
the divisor $D\subset Y$ is not $\Q$-Cartier. It is the {\em midpoint of
a Sarkisov link} (compare \cite{CPR}, 4.1 (3)); we develop this idea in
Part~II. The ideal case is when each $\Ga_i\subset X_1$ is a copy of $\PP^1$
with normal bundle $\Oh(-1,-1)$, or equivalently, $Y$ has only
ordinary nodes on $D$. We prove that this happens generically in all
our families.

In other situations, Type~I allows $\fie$ to be an isomorphism, typically
for $X$ of large index. At the other extreme, the Type~I condition on its
own does not imply that $-K_{X_1}$ is big, and $\fie$ could be an elliptic
Weierstrass fibration over $D=\PP(1,a,r-a)$, although this never happens
for codimension~4 Fano 3-folds. Also $\fie$ might
contract a surface to a curve of canonical singularities of $Y$; then
$X\broken Y$ is a ``bad link'' in the sense of \cite{CPR}, 5.5. We know
examples of this if $X$ is not required to be $\Q$-factorial and prime,
but none with these conditions.

\begin{exun} \rm
Consider the general codimension~2 complete intersection
\begin{equation}
X_{12,14}\subset\PP(1,1,4,6,7,8)_{\Span{x,a,b,c,d,e}}.
\end{equation}
The coordinate point $P_e=(0,\dots,0,1)$ is necessarily contained in
$X$: near it, the two equations $f_{12}:be=F_{12}$ and
$g_{14}:ce=G_{14}$ express $b$ and $c$ as implicit functions of the
other variables, so that $X$ is locally the orbifold point
$\frac18(1,1,7)$ with orbinates $x,a,d$.

Eliminating $e$ from $f_{12},g_{14}$ projects $X_{12,14}$ birationally
to the hypersurface
$Y_{18}:(bG-cF=0)\subset\PP(1,1,4,6,7)_{\Span{x,a,b,c,d}}$. Note that
$Y$ contains the plane $D=\PP(1,1,7)_{\Span{x,a,d}}=V(b,c)$, and has in
general $24=\frac17\times12\times14$ ordinary nodes at the points
$F=G=0$ of $D$.

In this case, the Kustin--Miller unprojection of the ``opposite''
divisor $(b=F=0)\subset Y$ completes the 2-ray game on $X_1$ to a
Sarkisov link, in the style of Corti and Mella \cite{CM}: the flop
$X_1\to Y\ot Y^+$ blows this up to a $\Q$-Cartier divisor, and
the unprojection variable $z_2=c/b=G/F$ then contracts it to a
nonorbifold terminal point $P_z\in
Z_{14}\subset\PP(1,1,4,7,2)_{\Span{x,a,b,d,z}}$.
\end{exun}

\subsection{Main theorem}

Write $P\in X$ for the numerical type of a codimension~4 Fano 3-fold of
index~1 marked with a Type~I centre. There are 115 or 116 candidates for
$X$ (depending on how you count the initial case); some have two or
three centres, and treating them separately makes 162 cases for $P\in X$.

\begin{thm}\label{th!main}
Let $P\in X$ be as above; then the projected variety is realised as a
codimension~$3$ Fano $Y\subset\wP^6$, and $Y$ can be made to contain a
co\-or\-dinate stratum $D=\PP(1,a,r-a)$ of $\wP^6$ in several ways.

For every numerical case $P\in X$, there are several formats, at least
one Tom and one Jerry (see Definition~\ref{d!TJ}) for which the general
$D\subset Y$ only has nodes
on $D$, and unprojects to a quasi\-smooth Fano $3$-fold $X\subset\wP^7$.
In different formats, the resulting $Y$ have different numbers of nodes
on $D$, so that the unprojected quasi\-smooth varieties $X$ have
different Betti numbers. Therefore in each of the $115$ numerical cases
for $X$, the Hilbert scheme has at least two components containing
quasismooth Fano $3$-folds.
\end{thm}

\subsection{Discussion of the result}
The theorem constructs around 320 different families of quasi\-smooth
Fano \hbox{$3$-folds}. We do not burden the journal pages with the
detailed lists, the {\em Big Table} in the Graded Ring Database
\cite{GRDB}; the case worked out in Section~\ref{exa!main} may be
adequate for most readers. Our data and the software tools for
manipulating them are available from \cite{GRDB}.

Our 162 cases for $P\in X$ project to $D\subset Y\subset w\PP^6$; of the
69 codimension~3 families of Fanos $Y$ that are $5\times5$ Pfaffians, 67
are the images of projections, each having up to four candidate planes
$D\subset Y$. For each of the 162 candidate pairs $D\subset Y$, we study
5 Tom and 10 Jerry formats, of which at least one Tom and one Jerry
succeeds (often one more, occasionally two), so that
Theorem~\ref{th!main} describes around 450 constructions of pairs $P\in
X$ of quasismooth Fano $3$-folds with marked centre of projection,
giving around 320 different families of $X$.

Theorem~\ref{th!main} covers codimension~4 Fano 3-folds of index~1 for
which there exists a Type~I centre. If one believes the possible
conjecture raised in \cite{ABR}, 4.8.3 that every Fano 3-fold in the
Mori category (that is, with terminal singularities) admits a
$\Q$-smoothing, this also establishes the components of the Hilbert
scheme of codimension~4 Fano 3-folds in these numerical cases. The main
novelty of this paper (and this was a big surprise to us) is that in
every case, the moduli space has 2, 3 or 4 different components.

An important remaining question is which $X$ are prime. In some cases,
our Tom or Jerry matrices have a zero entry, possibly after massaging.
Then 3 of the Pfaffian equations are binomial, which implies that $X$
has class group of rank $\rho\ge2$. This happens in the ancestral
examples of Section~\ref{sec!anc} and the easier cases
\ref{s!t2}--\ref{s!j25} of Section~\ref{exa!main}. Our Big Table
confirms that if we set aside all these cases with a zero, each of our
numerical possibilities for Type~I centres $P\in X$ admits exactly one
Tom and one Jerry construction that is potentially prime. Compare
Takagi's cases discussed in Section~\ref{s!intr}. We return to this
question in Part~II.

\subsection{Flowchart} \label{ss!flchart}
Our proof in Sections~\ref{s!fail}--\ref{s!grdb} applies computer algebra
calculations and verifications to a couple of thousand cases; any of
these could in principle be done by hand. We go to the database for
candidates for $P\in X$, figure out the weights of the coordinates of
$D\subset Y\subset w\PP^6$ and the matrix of weights, and list all
inequivalent Tom and Jerry formats. Section~\ref{s!fail} gives criteria
for a format to fail. In the cases that pass these tests,
Section~\ref{s!ns} contains an algorithm to produce $D\subset Y$ in the
given format, and to prove that it has only allowed singularities (that
is, only nodes on~$D$). Section~\ref{s!Jac} contains the Chern class
calculation for the number of nodes, so proving that the different
constructions build topologically distinct varieties.
Section~\ref{s!grdb} gives ``quick start-up'' instructions; do not
under any circumstances read the \verb!README! file.

\subsection{Further outlook}
The reducibility phenomenon appearing in this paper is characteristic
of Gorenstein in codimension~4; we have several current preprints and
work in progress addressing different aspects of this. See for example
\cite{Ki}.

This paper concentrates on 115 numerical cases of codimension~4 Fano
\hbox{3-folds} of index~1. Most of the remaining numerical cases from
Alt{\i}nok's list of 145 \cite{APhD} can be studied in terms of more
complicated Type~II or Type~IV unprojections, when the unprojection
divisor is not projectively normal; see \cite{Ki} for an introduction.
We believe that codimension~5 is basically similar: most cases have
two or more Type~I centres that one can project to smaller codimension,
leading to parallel unprojection constructions.

The methods of this paper apply also to other categories of varieties,
most obviously K3 surfaces and Calabi--Yau 3-folds. K3 surfaces are
included as general elephants $S\in|{-}K_X|$ in our Fano 3-folds,
although the K3 is unobstructed, so that passing to the elephant hides
the distinction between Tom and Jerry. We can also treat some of the
Fano 3-folds of index $>1$ of Suzuki's thesis \cite{S}, \cite{BS}; we
have partial results on the existence of some of these families, and
hope eventually to cover the cases not excluded by Prokhorov's
birational methods \cite{Pr}.

This paper uses Type~I projections $X\broken Y$ to study the biregular
question of the existence and moduli of $X$; however, in each case, the
Kawamata blowup $X_1\to X$ initiates a 2-ray game on $X_1$, with the
anticanonical model $X_1\to Y$ and its flop $Y\ot Y^+$ as first step. In
many cases, we know how to complete this to a Sarkisov link using Cox
rings, in the spirit of \cite{CPR}, \cite{CM}, \cite{BCZ} and \cite{BZ};
we return to this in Part~II.

\section{Extended example} \label{exa!main}
The case $g=0$ plus basket
$\bigl\{\frac12(1,1,1),\frac13(1,1,2),\frac14(1,1,3),\frac15(1,1,4)\bigr
\}$ gives the codimension~4 candidate $X\subset\PP^7(1,1,2,3,3,4,4,5)$
with Hilbert numerator
\begin{equation}
1-2t^6-3t^7-3t^8-t^9+t^9+4t^{10}+6t^{11}+\cdots+t^{22}.
\end{equation}
It has three different possible Type~I centres, namely the $\frac13$,
$\frac14$ or $\frac15$ points. We project away from each of these,
obtaining consistent results; each case leads to four unprojection
constructions for $X$, two Toms and two Jerries:
\begin{description}

\item{from $\frac13$:} gives $\PP(1,1,2)\subset Y\subset\PP(1,1,2,3,4,4,5)$
with matrix of weights
\begin{equation}
\begin{pmatrix}
2&2&3&4 \\
&3&4&5 \\
&&4&5 \\
&&&6
\end{pmatrix} \qquad\hbox{and}\qquad
\begin{array}l
\hbox{Tom$_2$ has 13 nodes} \\
\hbox{Tom$_1$ has 14 nodes} \\
\hbox{Jer$_{45}$ has 16 nodes} \\
\hbox{Jer$_{25}$ has 17 nodes}
\end{array} \label{eq!3.2}
\end{equation}

\item{from $\frac14$:} gives $\PP(1,1,3)\subset
Y\subset\PP(1,1,2,3,3,4,5)$ with matrix of weights
\begin{equation}
\begin{pmatrix}
2&3&3&4 \\
&3&3&4 \\
&&4&5 \\
 &&&5
\end{pmatrix} \qquad\hbox{and}\qquad
\begin{array}l
\hbox{Tom$_3$ has 9 nodes} \\
\hbox{Tom$_1$ has 10 nodes} \\
\hbox{Jer$_{35}$ has 12 nodes} \\
\hbox{Jer$_{15}$ has 13 nodes}
\end{array} \label{eq!3.3}
\end{equation}

\item{from $\frac15$:} gives $\PP(1,1,4)\subset
Y\subset\PP(1,1,2,3,3,4,4)$ with matrix of weights
\begin{equation}
\begin{pmatrix}
2&2&3&3 \\
&3&4&4 \\
&&4&4 \\
 &&&5
\end{pmatrix} \qquad\hbox{and}\qquad
\begin{array}l
\hbox{Tom$_4$ has 8 nodes} \\
\hbox{Tom$_2$ has 9 nodes} \\
\hbox{Jer$_{24}$ has 11 nodes} \\
\hbox{Jer$_{14}$ has 12 nodes}
\end{array} \label{eq!3.4}
\end{equation}
\end{description}
Specifically, we assert that {\em in each of these 12 cases, if we pour
general elements of the ideal $I_D$ and general elements of the ambient
ring into the Tom or Jerry matrix $M$ as specified in
Definition~\ref{d!TJ}, the Pfaffians of $M$ define a Fano $3$-fold $Y$
having only the stated number of nodes on $D$, and the resulting $X$ is
quasismooth.} Section~\ref{s!ns} verifies this claim by cheap computer
algebra, although we work out particular cases here without such
assistance. Section~\ref{s!Jac} computes the number of nodes in each
case from the numerical data. Imposing the unprojection plane $D$ on the
general quasismooth $Y_t$ introduces singularities on $Y=Y_0$, nodes in
general, which are then resolved on the quasismooth $X_1$. Each node
thus gives a conifold transition, replacing a vanishing cycle $S^3$ by a
flopping line $\PP^1$, and therefore adds 2 to the Euler number of $X$;
so the four different $X$ have different topology.

The unprojection formats and nonsingularity algorithms establish the
existence of four different families of quasismooth Fano 3-folds $X$.
The rest of this section analyses these in reasonably natural formats;
an ideal would be to free ourselves from unprojection and computer
algebra, although we do not succeed completely.

For illustration, work from $\frac13$; take
$X\subset\PP^7(1,1,2,3,3,4,4,5)_{\Span{x,a,b,c,d,e,f,g}}$, and assume
that $P_d=(0,0,0,0,1,0,0,0)$ is a Type~1 centre on $X$ of type
$\frac13(1,1,2)$. The assumption means that $P\in X$ is quasismooth with
orbinates $x,a,b$. The cone over $X$ is thus a manifold along the
$d$-axis, and therefore, by the implicit function theorem, four of the
generators of $I_X$ form a regular sequence locally at $P_d$, with
independent derivatives, say $cd=\cdots$, $de=\cdots$, $df=\cdots$,
$dg=\cdots$ of degrees $6,7,7,8$. Eliminating $d$ gives the Type~I
projection $X\broken Y$ where $Y\subset\PP^6(1,1,2,3,4,4,5)$ has Hilbert
numerator
\begin{equation}
1-t^6-t^7-2t^8-t^9+t^{10}+2t^{11}+t^{12}+t^{13}-t^{19}.
\end{equation}
Let $Y$ be a $5\times5$ Pfaffian matrix with weights as in
\eqref{eq!3.2}. Since rows 2 and 3 have the same weights, we can
interchange the indices 2 and 3 throughout; thus Tom$_2$ is equivalent
to Tom$_3$, Jer$_{25}$ to Jer$_{35}$, and so on.

\subsection{Failure} \label{s!fa}
Some Tom and Jerry cases fail, either for coarse or for more subtle
reasons; for example, it sometimes happens that for reasons of weight,
one of the variables $x_i$ cannot appear in the matrix, so the variety
is a cone, which we reject. Section~\ref{s!fail} discusses failure
systematically.

In the present case $D=\PP(1,1,2)_{\Span{x,a,b}}$, the generators of
$I_D=(c,e,f,g)$ all have weight $\ge3$, but $\wt m_{12},m_{13}=2$. Thus
requiring $m_{12},m_{13}\in I_D$ forces them to be zero, making the
Pfaffians $\Pf_{12.34}$ and $\Pf_{12.35}$ reducible. This kills Tom$_4$,
Tom$_5$, Jer$_{1i}$ for any $i$ and Jer$_{23}$. The same argument
says that Tom$_2$ has $m_{13}=0$ and Jer$_{25}$ has $m_{12}=0$, a key
simplification in treating them: a zero in $M$ makes three of the
Pfaffians binomial.

We see below that Jer$_{24}$ fails for an interesting new reason. The
other cases all work, as we could see from the nonsingularity algorithm
of Section~\ref{s!ns}. Tom$_2$ and Jer$_{25}$ are simpler, and we start
with them, whereas Tom$_1$ and Jer$_{45}$ involve heavier calculations;
they are more representative of constructions that possibly lead to
prime $X$.

\subsection{Tom$_2$}\label{s!t2} The analysis of the matrix proceeds as:
\begin{equation}
\renewcommand{\arraycolsep}{.35em}
\begin{pmatrix}
K_2 & 0 & c & e \\
& L_3 & M_4 & N_5 \\
&& f & g \\
&&& \Span{c,e,f,g}_6
\end{pmatrix} \mapsto
\begin{pmatrix}
b & 0 & c & e \\
& L_3 & M_4 & N_5 \\
&& f & g \\
&&& 0
\end{pmatrix}
\mapsto
\begin{pmatrix}
b & c & e \\
d & M & N \\
L & f & g
\end{pmatrix}
\label{eq!2.8}
\end{equation}
here $m_{13}=0$ is forced by low degree, $K_2$, $L_3$, $M_4$, $N_5$ are
general forms of the given degrees, that we can treat as tokens
(independent indeterminates), and the four entries
$m_{14},m_{15},m_{34},m_{35}$ are general elements of $I_D$ that we
write $c,e,f,g$ by choice of coordinates. Next, $m_{45}$ can be whittled
away to 0 by successive row-column operations that do not harm the
remaining format; seeing this is a ``crossword puzzle'' exercise that
uses the fact that $m_{13}=0$ and all the entries in Row~2 are general
forms. For example, subtracting a suitable multiple of Row~1 from Row~5
(and then the same for the columns) kills the $c$ in $m_{45}$, while
leaving $m_{15}$ and $m_{35}$ unchanged (because $m_{11}=m_{13}=0$) and
modifying $N_5$ by a multiple of $K_2$, which is harmless because $N_5$
is just a general ring element of weight~5.

The two zeros imply that all the Pfaffians are binomial, and, as in
\ref{s!same}, putting in the unprojection variable $d$ of weight 4
gives the $2\times2$ minors of the matrix on the right. The equations
describe $X$ inside the projective cone over
$w(\PP^2\times\PP^2)\subset\PP(2,3^3,4^3,5^2)$ with vertex
$\PP^1_{\Span{x,a}}$ as the complete intersection of three general forms
of degree $3,4,5$ expressing $L,M,N$ in terms of the other variables.
(It is still considerably easier to do the nonsingularity computation
after projecting to smaller codimension.)

\subsection{Jer$_{25}$}\label{s!j25} We start from
\begin{equation}
\begin{pmatrix}
0 & b & L_3 & f \\
& c & e & g \\
&& M_4 & \la_1 e \\
&&& \mu_3 c + \nu_2 e
\end{pmatrix}
\label{eq!2.9}
\end{equation}
where $m_{12}=0$ is forced by low degree, and we put tokens $b,L,M$ in
place of the free entries $m_{13},m_{14},m_{34}$. We have cleaned out
$m_{35}$ and $m_{45}$ as much as we can; the quantities
$b,L,M,\la,\mu,\nu$ are general ring elements of the given weights.

We have to adjoin $d$ together with unprojection equations for
$dc,de,df,dg$. There are various ways of doing this, including the
systematic method of writing out the Kustin--Miller homomorphism between
resolution complexes, that we use only as a last resort. An ad hoc
parallel unprojection method is to note that $g$ appears only as the
entry $m_{25}$, so we can project it out to a codimension~2 c.i.\
containing the plane $c=e=f=0$:
\begin{equation}
\begin{pmatrix} \mu b & \nu b-\la L & M \\ L & -b & 0 \end{pmatrix}
\begin{pmatrix} c \\ e \\ f \end{pmatrix}=0.
\end{equation}
The equations for $dc,de,df$ come from Cramer's rule, and we can
write the unprojection in rolling factors format:
\begin{equation}
\bigwedge^2\begin{pmatrix}
b & L & f & d\\
c & e & g & M
\end{pmatrix}
\quad\hbox{and}\quad
\begin{array}l
\mu b^2 + \nu bL-\la L^2+ df,\\
\mu bc + \nu cL-\la eL + Mf, \\
\mu c^2 + \nu ce-\la e^2 + Mg.
\end{array}
\label{eq!rf}
\end{equation}
The first set of equations of \eqref{eq!rf}, with the entries viewed as
indeterminates, defines
$w(\PP^1\times\PP^3)\subset\PP(2,3,3,3,4,4,4,5)_{\Span{b,c,d,L,e,f,M,g}}$;
the second set is a single quadratic form evaluated on the rows, so
defines a divisor in the cone over this with vertex
$\PP^1_{\Span{x,a}}$. Finally, setting $L,M$ general forms gives $X$ as
a complete intersection in this.

\subsection{Jer$_{24}$ fails}
The matrix has the form
\begin{equation}
\begin{pmatrix}
0 & b & c & L_4 \\
& c & f & g \\
&& e & M_5 \\
&&& \Span{c,e,f,g}_6
\end{pmatrix}
\mapsto
\begin{pmatrix}
0 & b & c & L_4 \\
& c & f & g \\
&& e & M_5 \\
&&& 0
\end{pmatrix}
\end{equation}
The entries in the rows and columns through the pivot $m_{24}=f$ are
general elements of the ideal $I_D=(c,e,f,g)$. As before, $m_{12}=0$ is
forced by degrees. Although \ref{s!fish}, (5) fails this for a
mechanical reason, we discuss it in more detail as an instructive case,
giving a perfectly nice construction of the unprojected variety $X$,
that happens to be slightly too singular. First, please check that the
entry $m_{45}$ can be completely taken out by row and column operations.
For example, to get rid of the $e$ term in $m_{45}$, add $\al_3$ times
Row~3 to Row~5; in $m_{25}$ this changes $g$ to $g+\al c$, that we
rename $g$.

One sees that the equations of the unprojected variety $X$ take the form
\begin{equation}
\bigwedge^2\begin{pmatrix} b&c&e&f \\ d&L&M&g \end{pmatrix}=0
\quad\hbox{and}\quad \left\{
\begin{matrix}
bf=c^2, \\ bg=cL, \\ dg=L^2.
\end{matrix}
\right.
\end{equation}
(exercise, hint: project out $f$ or $g$). In straight projective
space, these equations define $\PP^1\times Q\subset\PP^1\times\PP^3$
where $Q\subset\PP^3$ is the quadric cone. This is singular in
codimension~2, so the 3-fold $X$ cannot have isolated singularities.

\subsection{Tom$_1$}
The matrix and its clean form are
\begin{equation}
\begin{pmatrix}
b& K_2 & L_3 & M_4 \\
& c & e & g \\
&& f & \Span{c,e,f,g}_5 \\
&&& \Span{c,e,f,g}_6
\end{pmatrix}
\mapsto
\begin{pmatrix}
b& K & L & M \\
& c & e & g \\
&& f & \la_1e \\
&&& \mu_3c+\nu_2e
\end{pmatrix}
\end{equation}
where $K,L,M$ and $\la,\mu,\nu$ are general forms, that we treat as
tokens. We add a multiple of Column~2 to Column~5 to clear $c$ from
$m_{35}$, so we cannot use the same operation to clear $e$ from
$m_{45}$. The nonsingularity algorithm of Section~\ref{s!ns} ensures
that for general choices this has only nodes on $D$.

We show how to exhibit $X$ as a triple parallel unprojection from a
hypersurface in the product of three codimension~2 c.i.\ ideals (compare
\ref{s!xtra}). Since $g$ only appears as $m_{25}$, it is
eliminated by writing the two Pfaffians $\Pf_{12.34}$ and $\Pf_{13.45}$
as:
\begin{equation}
\begin{pmatrix} L & -K & b \\ \mu K & \nu K-\la L & M \end{pmatrix}
\begin{pmatrix} c \\ e \\ f \end{pmatrix}= 0;
\end{equation}
in the same way, $\Pf_{12.45}$ and $\Pf_{12.35}$ eliminate $f$:
\begin{equation}
\begin{pmatrix} M & \la b & -K \\ \mu b & M+\nu b &-L \end{pmatrix}
\begin{pmatrix} c \\ e \\ g \end{pmatrix} = 0.
\end{equation}
Cramer's rule applied to these gives the unprojection equations for $d$:
\begin{equation}
\begin{array}{l}
dc = KM + \nu bK -\la bL, \\
de = LM - \mu bK,
\end{array} \quad
\begin{array}{l}
df = -\mu K^2 + \nu KL - \la L^2, \\
dg = M^2 + \nu bM - \la\mu b^2.
\end{array}
\end{equation}
The combination eliminating $d$, $f$ and $g$ is
\begin{equation}
eKM - cLM - \la beL + \mu bcK + \nu beK = 0.
\label{eq!Z}
\end{equation}
This is a hypersurface $Z_{10} \subset
\PP^4(1,1,2,3,4)_{\Span{x,a,b,c,e}}$ contained in the product ideal of
$I_d=(c,e)$, $I_f=(b,M_4)$, $I_g=(K_2,L_3)$. The unprojection planes
$\Pi_d$, $\Pi_f$, $\Pi_g$ are projectively equivalent to $\PP(1,1,2)$,
$\PP(1,1,3)$, $\PP(1,1,4)$, but we cannot normalise all three of them to
coordinate planes at the same time. Their pairwise intersection is:
\begin{align*}
\Pi_d\cap \Pi_f &= \hbox{the 4 zeros of $M_4$ on the line $b=c=e=0$,} \\
\Pi_d\cap \Pi_g &= \hbox{the 3 zeros of $L_3$ on the line $c=e=K=0$,} \\
\Pi_f\cap \Pi_g &= \hbox{the 2 zeros of $K_2$ on the line $b=L=M=0$.}
\end{align*}

\paragraph{Nonsingularity based on \eqref{eq!Z}}
All the assertions we need for $Y$ and $X$ are most simply derived from
\eqref{eq!Z}. The linear system $|I_d\cdot I_f\cdot
I_g\cdot\Oh_{\PP}(10)|$ of hypersurfaces through the three unprojection
planes has base locus the planes themselves, together with the curve
$(b=c=K_2=0)$, which is in the base locus because the term $eLM\in
I_d\cdot I_f\cdot I_g$ has degree~11 and so does not appear in the
equation of $Z$. This curve is a pair of generating lines
$(K=0)\subset\PP(1,1,4)_{\Span{x,a,e}}$. One sees that for general
choices, one of the terms $cLM$ or $\la beL$ in $Z$ provides a nonzero
derivative $LM$ or $\la eL$ at every point along this curve away from
the three planes.

The singular locus of $Z$ on $\Pi_d=\PP(1,1,2)$ is given by
\begin{equation}
\frac{\partial Z}{\partial c}=-LM+\mu bK =0, \quad
\frac{\partial Z}{\partial e}=KM-\la bL+\nu bK=0.
\end{equation}
For general choices, these are $21=\frac{7\times6}2$ reduced points of
$\PP(1,1,2)$, including the 4 points of $\Pi_d\cap \Pi_f$ and the 3
points of $\Pi_d\cap \Pi_g$; after unprojecting $\Pi_f$ and $\Pi_g$,
this leaves 14 nodes of Tom$_1$, as we asserted in \eqref{eq!3.2}. The
calculations for the other planes are similar.

We believe that $Z_{10} \subset \PP^4(1,1,2,3,4)$ has class group
$\Z^4$ generated by the hyperplane section $A=-K_Z$ and the three planes
$\Pi_d$, $\Pi_f$, $\Pi_g$, so that $X$ is prime.

\subsection{Jer$_{45}$} \label{s!J45}
The tidied up matrix is
\begin{equation}
\begin{pmatrix}
b & -L_2 & c & e \\
& M_3 & e & g \\
&& f & \la_2 c \\
 &&& m_{45}
\end{pmatrix},
\label{eq!3.12}
\end{equation}
with pivot $m_{45}=\de_3c+\ga_2e+\be_2f+\al_1g$; we use row and column
operations and changes of coordinates in $I_D=(c,e,f,g)$ to clean $c$
and $f$ out of $m_{24}$, but we cannot modify the pivot $m_{45}$ without
introducing multiples of $b,L,M$ into Row~4 or Row~5, spoiling the
Jer$_{45}$ format.

We get parallel unprojection constructions for $X$ by eliminating $f$ or
$g$ or both. First, subtract $\al$ times Row~2 from Row~4, and ditto
with the columns, to take $g$ out of $m_{45}$. This spoils the format by
$c\mapsto c-\al b\notin I_D$ in $m_{14}$, but does not change the
Pfaffian ideal. The new matrix only contains $g$ in $m_{25}$;
the two Pfaffians not involving it are $\Pf_{12.34}$ and the modified
$\Pf_{13.45}$, giving
\begin{equation}
\begin{pmatrix} M & L & b \\
\de L+\la c-\al\la b & \ga L-\al M & \be L-e \end{pmatrix}
\begin{pmatrix} c \\ e \\ f \end{pmatrix}=0.
\label{eq!gless}
\end{equation}
Eliminating $f=m_{34}$ is similar, with $\Pf_{12.35}$ and modified
$\Pf_{12.45}$ giving
\begin{equation}
\begin{pmatrix} \la b & M & L \\
\de b-\be M & \ga b+e-\be L & \al b-c \end{pmatrix}
\begin{pmatrix} c \\ e \\ g \end{pmatrix}=0.
\label{eq!fless}
\end{equation}
We derive the unprojection equations for $d$ using Cramer's rule:
\begin{equation}
\begin{array}l
dc=-L(e-\be L)-\ga Lb+\al Mb, \\
de=M(e-\be L)+\la b(c-\al b)+\de Lb, \\
df=-\la L(c-\al b)-\de L^2+\ga LM-\al M^2, \\
dg=\la b(e-\be L)+M(g-\de b)+\ga\la b^2+\be M^2.
\end{array}
\label{eq!deqns}
\end{equation}

This is also a triple parallel unprojection, but with a difference: the
hypersurface $Z_{10}\subset\PP(1,1,2,3,4)$ obtained by eliminating $f$
from \eqref{eq!gless} or $g$ from \eqref{eq!fless} or $d$ from the first
two rows of \eqref{eq!deqns} is now
\begin{equation}
e(e-\be L)L+\de cbL+\ga ebL+\la bc(c-\al b)+M(ce-\be cL-\al be)=0.
\label{eq!H}
\end{equation}
It is in the inter\-section of the three codimension~2 c.i.\
unprojection ideals $I_d=(c,e)$, $I_f=(b,e-\be L)$, $I_g=(c-\al b,L)$,
but not in their product: the first 4 terms are clearly in the product
ideal. The interesting part is the bracket in the last term, which
cannot be in the product since it has terms of degree~2, but is in
$I_d\cap I_f\cap I_g$, because
\begin{equation}
c(e-\be L)-\al be = e(c-\al b)-\be Lc.
\end{equation}
The slogan is {\em like lines on a quadric}; the three ideals have
linear combinations of $b,c$ as first generator, and of $e,L$ as second
generator, like three disjoint lines $x=z=0$, $y=t=0$ and $x=t,y=z$ on
$Q:(xy=zt)$. One analyses the singularities of $Z_{10}$ from this much
as before; we believe that $\Cl Z=\Span{A,D_1,D_2,D_3}$, so that the
triple unprojection $X$ is prime.

\section{Failure}\label{s!fail}

We give reasons for failure following the introductory discussion in
Section~\ref{exa!main}; we don't need to treat all the possible tests in
rigorous detail, or the logical relations between them. For the
structure of our proof, the point of this section is merely to give
cheap preliminary tests to exclude all the candidates $D\subset Y$ that
will not pass the nonsingularity algorithm in Section~\ref{s!ns}.

\subsection{Easy fail at a coordinate point} \label{s!easy}

Consider a coordinate point $P_i=P_{x_i}\in Y$. In either of the
following cases, $P_i$ cannot be a hyperquotient point, let alone
terminal, and we can safely fail the candidate $D\subset Y$:

\begin{enumerate}
\renewcommand{\labelenumi}{(\arabic{enumi})}

\item $x_i$ does not appear in the matrix $M$.

\item $x_i$ does not appear as a pure power in any entry of $M$, which
thus has rank zero at $P_i$.

\end{enumerate}

\subsection{Fishy zero in $M$ and excess singularity} \label{s!fish}

Suppose we can arrange that $m_{12}=0$, if necessary after row and
column operations; then the subscheme $Z=V(\{m_{1i},m_{2i}\mid
i=3,4,5\})$ is in the singular locus of $Y$. Indeed, the three Pfaffians
$\Pf_{12,ij}$ are in $I_Z^2$, so do not contribute to the Jacobian at
points of $Z$. The case that $\dim Z=0$ and $Z\subset D$ is perfectly
acceptable and happens in a fraction of our successful constructions
(see Tom$_2$ and Jer$_{25}$ in Section~\ref{exa!main}). Notice that
$\dim Z=0$ if and only if the 6 forms $m_{1i},m_{2i}$ make up a regular
sequence for $\PP^6$; in the contrary case, the zero is {\em fishy}.
Thus any little coincidence between the six $m_{1i},m_{2i}$ fails
$D\subset Y$. The tests we implement are:

\begin{enumerate}
\renewcommand{\labelenumi}{(\arabic{enumi})}
\setcounter{enumi}{2}
\item Two collinear zeros in $M$; see \ref{s!fa} for an example.

\item Two of the $m_{1i},m_{2i}$ coincide; see Section~\ref{exa!main},
Jer$_{24}$.

\item An entry $m_{1i}$ or $m_{2i}$ is in the ideal generated by the
other five.

\end{enumerate}

In fact, the tricky point here is how to read our opening ``Suppose we
can arrange that $m_{12}=0$''. The row and column operations clearly need
a modicum of care to preserve the format (i.e., the entries we require
to be in $I_D$). The harder point is that we may need a particular
change of basis in $I_D$ for the zero to appear. For example, in the
Tom$_5$ format for $\PP^2\subset Y\subset\PP(1^6,2)$, with matrix of
weights
$\begin{smallmatrix}1&1&1&2\\&1&1&2\\&&1&2\\&&&2\end{smallmatrix}$, the
lowest degree Pfaffian is quadratic in three variables of weight~1, so
we can write it $xy-z^2$. Mounting this as a Pfaffian in these
coordinates, we can force a fishy zero, with two equal entries $z$
arising from the term $z^2$. (The same applies to several candidates,
but this is the only one that fails solely for this reason.)

\subsection{More sophisticated and ad hoc reasons for failure}

For the unprojected $X$ to have terminal singularities, $Y$ itself must
also: it is the anticanonical model of the weak Fano 3-fold $X_1$. We
can test for this at a coordinate point $P$ of index $r>1$: by Mori's
classification, $Y$ is either quasismooth at $P$, or a hyperquotient
singularity with local weights $\frac1r(1,a,r-a,0)$ or
$\frac14(1,1,3,2)$. Thus we can fail the candidate $D\subset Y$ if:
\begin{enumerate}
\renewcommand{\labelenumi}{(\arabic{enumi})}
\setcounter{enumi}{5}
\item A coordinate point off $D$ is a nonterminal hyperquotient singularity.

\item A coordinate point on $D$ is a nonterminal hyperquotient singularity.

\end{enumerate}
These tests dispatch most of the remaining failing candidates.

\begin{enumerate}
\renewcommand{\labelenumi}{(\arabic{enumi})}
\setcounter{enumi}{7}
\item Ad hoc fail. Just two cases have nonisolated singularities not
revealed by the elementary tests so far:
\begin{enumerate}

\item Tom$_4$ for $\PP(1,2,3)\subset Y\subset\PP(1^2,2,3^2,4^2)$ with weights
$\begin{smallmatrix}2&2&3&3\\&3&4&4\\&&4&4\\&&&5\end{smallmatrix}$\,;

\item Jer$_{12}$ for $\PP(1,2,3)\subset Y\subset\PP(1^2,2^2,3^2,4)$ with weights
$\begin{smallmatrix}2&2&2&3\\&3&3&4\\&&3&4\\&&&4\end{smallmatrix}$\,.
\end{enumerate}
Each of these has a $\frac{1}2(1,1,1,0;0)$ hyperquotient singularity at
the $\frac12$ point of $D$. Such a point may be terminal if it is an
isolated double point, but the format of the matrix prevents this. The
second case also fails at the index~4 point $P_7$ lying off $D$: it is a
hyperquotient singularity of the exceptional type $\frac{1}4(1,1,3,2;2)$
with the right quadratic part to be terminal. However, it lies on a
curve of double points along the line $\PP(2,4)$ joining $P_7$ to the
$\frac12$ point on $D$: in local coordinates $x,a,e,b$ at $P_7$, the
equation is $xa=e^2+b\times\hbox{terms in $(x,a,e)^2$}$.

\end{enumerate}

\section{Nonsingularity and proof of Theorem~\ref{th!main}} \label{s!ns}

To prove Theorem~\ref{th!main}, we need to run through a long list of
candidate \hbox{3-folds} $D\subset Y \subset w\PP^6$ with choice of
format Tom$_i$ or Jer$_{ij}$. We exclude many of these by the automatic
methods of Section~\ref{s!fail}. In every remaining case, we run a
nonsingularity algorithm to confirm that the candidate can be
unprojected to a codimension~4 Fano 3-fold $X$ with terminal
singularities (in fact, we conclude also quasismooth). For the proof of
Theorem~\ref{th!main}, we check that at least one Tom and one Jerry
works for each case $D\subset Y$.

We outline the proof as a pseudocode algorithm; our implementation is
discussed in Section~\ref{s!grdb}. The justification of the algorithm is
that it works in practice. A priori, it could fail, e.g., the singular
locus of $Y$ on $D$ could be more complicated than a finite set of
nodes, or all three coordinate lines of $D$ could contain a node, but by
good luck such accidents never happen.

\subsection{Nonsingularity analysis} \label{s!sings}

We work with any $D\subset Y$ not failed in Section~\ref{s!fail}. The
homogeneous ideal $I_Y$ is generated by the $4\times4$ Pfaffians of $M$.
Differentiating the 5 equations $\Pf$ with respect to the seven
variables gives the $5\times7$ Jacobian matrix $J(\Pf)$. Its ideal
$I_{\Sing Y}=\bigwedge^3J(\Pf)$ of $3\times3$ minors defines the
singular locus of $Y$; more precisely, it generates the ideal sheaf
$\sI_{\Sing Y}\subset\Oh_{\PP^6}$. Our claim is that the only
singularities of $Y$ lie on $D$, and are nodes. For this, we check that
\begin{enumerate}
\renewcommand{\labelenumi}{(\alph{enumi})}
\item $\Sing Y\subset D$, or equivalently $I_D\subset\Rad(I_{\Sing Y})$.

\item The restriction $\sI_{\Sing Y}\cdot\Oh_D$ defines a reduced
subscheme of $D$.

\end{enumerate}
In fact (b) together with Lemma~\ref{l!nodes} imply that $Y$ has only
nodes. In practice, we may work on a standard affine piece of $D$
containing all the singular points: it turns out in every case that some
1-strata of $D$ is disjoint from the singular locus.

\subsection{Proof of Theorem~\ref{th!main}}

We start with the data for a candidate $P\in X\subset w\PP^7$: a genus
$g\ge-2$ and a basket $\sB$ of terminal quotient singularities, or
equivalently, the resulting Hilbert series (see \cite{ABR}). We give a
choice of 8 ambient weights $W_X$ of $w\PP^7$ and a choice of Type~I
centre $P=\frac1r(1,a,r-a)$ from the basket. The Type~I definition
predicts that the ambient weights of $Y\subset w\PP^6$ are $W_X
\setminus \{r\}$ and that $D=\PP(1,a,r-a)$ can be chosen to be a
coordinate stratum of $w\PP^6$. We analyse all possible Tom and Jerry
formats for $D\subset Y\subset w\PP^6$.

\paragraph{Step 1} Set up coordinates $x_1,x_2,x_3,x_4$, $y_1,y_2,y_3$
on $w\PP^6$; here $x_{1\dots4}$ is a regular sequence generating $I_D$,
and $y_1,y_2,y_3$ are coordinates on $D$.

\paragraph{Step 2} The numerics of \cite{CR} determine the weights
$d_{ij}$ of the $5\times5$ skew matrix $M$ from the Hilbert numerator of
$Y\subset w\PP^6$.

\paragraph{Step 3} Set each entry $m_{ij}$ of $M$ equal to a general
form, respectively a general element of the ideal $I_D$ of the given
degree $d_{ij}$, according to the chosen Tom or Jerry format (see
Definition~\ref{d!TJ}).

Tidy up the matrix $M$ as much as possible while preserving its Tom or
Jerry format. Some entries of $M$ may already be zero. Use coordinate
changes on $w\PP^6$ to set some entries of $M$ equal to single
variables. If possible, use row and column operations to simplify $M$
further. Check every zero of $M$ for failure for the mechanical reasons
discussed in \ref{s!fish}, followed by the other failing conditions of
\ref{s!easy}. Now any candidate that passes these tests actually works.

\paragraph{Step 4} Carry out the singularity analysis of \ref{s!sings}.

\paragraph{Step 5} Calculate the number of nodes as in
Section~\ref{s!Jac}; check that no two sets of unprojection data give
the same number of nodes.

\paragraph{Step 6 (optional)} Apply the Kustin--Miller algorithm
\cite{KM} to construct the equations of $X$. This is not essential to
prove that $X$ exists, but knowing the full set of equations is useful
if we want to put the equations in a codimension~4 format, for example
by projecting from another Type~I centre.

\section{Number of nodes} \label{s!Jac}

The unprojection divisor $D=V(x_{1\dots4})\subset \PP^6$ is a
codimension~4 c.i., with conormal bundle $\sI_D/\sI_D^2$ the direct sum
of four orbifold line bundles $\Oh_D(-x_i)$ on $D$. The ideal sheaf
$\sI_Y$ is generated by 5 Pfaffians that vanish on $D$, so each is
$\Pf_i=\sum a_{ij}x_j$. Thus the Jacobian matrix $\Jac$ restricted to
$D$ is the $5\times4$ matrix $(\bara_{ij})$, where bar is restriction
mod $I_D=(x_{1\dots4})$; the induced homomorphism to the conormal bundle
\begin{equation}
\sJ\colon \bigoplus_5\Oh_{\PP}(-\Pf_i)\onto \sI_Y/(\sI_D\cdot\sI_Y) \to
\sI_D/\sI_D^2
\label{eq!N}
\end{equation}
has generic rank~3. Its cokernel $\sN$ is the conormal sheaf to $D$ in
$Y$. It is a rank~1 torsion free sheaf on $D$ whose second Chern class
$c_2(\sN)$ counts the nodes of $Y$ on $D$. The more precise result is as
follows:

\begin{lem} \label{l!nodes}
\begin{enumerate}
\renewcommand{\labelenumi}{(\Roman{enumi})}

\item The cokernel $\sN$ is an orbifold line bundle at points of $D$
where $\rank\sJ=3$, that is, at quasi\-smooth points of $Y$.

\item Assume that $P\in D$ is a nonsingular point (not orbifold), and
that $\rank\sJ=2$ at $P$ and $=3$ in a punctured neighbourhood of $P$ in
$D$; then $\sN$ is isomorphic to a codimension~$2$ c.i.\ ideal $(f,g)$
locally at $P$. This coincides locally with the ideal
$\bigwedge^3\Jac\cdot\Oh_D$ generated by the $3\times3$ minors of the
Jacobian matrix.

\item Assume that $\bigwedge^3\Jac\cdot\Oh_D$ is reduced (locally the
maximal ideal $m_P$ at each point). Then $Y$ has an ordinary node at $P$.

\item If this holds everywhere then $c_2(\sN)$ is the number of
nodes of $Y$ on $D$.
\end{enumerate}
\end{lem}

\paragraph{Proof} The statement is the hard part; the proof is just
commutative algebra over a regular local ring. The rank~1 sheaf $\sN$ is
the quotient of a rank~4 locally free sheaf by the image of the
$5\times4$ matrix $\Jac=(\bara_{ij})$, of generic rank~3. It is a line
bundle where the rank is 3, and where it drops to~2, we can use a
$2\times2$ nonsingular block to take out a rank~2 locally free summand.
The cokernel is therefore locally generated by 2 elements, so is locally
isomorphic to an ideal sheaf $(f,g)$, a c.i.\ because the rank drops
only at $P$.

The minimal free resolution of $\sN$ is the Koszul complex of $f,g$; now
\eqref{eq!N} is also part of a free resolution of $\sN$, so covers the
Koszul complex. This means that the matrix $\Jac=(\bara_{ij})$ can be
written as its $2\times2$ nonsingular block and a complementary
$2\times3$ block of rank~1, whose two rows are $g\cdot v$ and $-f\cdot
v$ for $v$ a 3-vector with entries generating the unit ideal. Therefore
$\bigwedge^3\Jac$ generates the same ideal $(f,g)$.

If $(f,g)=(y_1,y_2)$ is the maximal ideal at $P\in D$ then the shape of
$\bigwedge^3\Jac$ says that two of the Pfaffians $\Pf_1,\Pf_2$ express
two of the variable $x_1,x_2$ as implicit functions; then a linear
combination $p$ of the remaining three has $\partial p/\partial x_3=y_1$
and $\partial p/\partial x_4=y_2$, so that $Y$ is a hypersurface with an
ordinary node at $P$. \QED\par\medskip

We now show how to resolve $\sN$ by an exact sequence involving direct
sums of orbifold line bundles on $D$, and deduce a formula for
$c_2(\sN)$.

\paragraph{Tom$_1$}
The matrix is
\begin{equation}
M = \begin{pmatrix}
K & L & M & N \\
& m_{23} & m_{24} & m_{25} \\
&& m_{34} & m_{35} \\
&&& m_{45}
\end{pmatrix}
\end{equation}
where $m_{ij}$ are linear forms in $x_{1\dots4}\in\sI_D$ with
coefficients in the ambient ring. When we write out $\Jac=(\bara_{ij})$,
the only terms that contribute are the derivatives $\partial/ \partial
x_{1\dots4}$, with the $x_i$ set to zero; thus only the terms that are
exactly linear in the $x_i$ contribute. Since $\Pf_1$ is of order $\ge2$
in the $x_i$, the corresponding row of the matrix $J$ is zero and we
omit it in \eqref{eq!J4}. Moreover, the first row $K,L,M,N$ of $M$
provides a syzygy $\Si_1=K\Pf_2+L\Pf_3+M\Pf_4+N\Pf_5\equiv0$ between the
4 remaining Pfaffians. Hence we can replace $J$ by the resolution
\begin{equation}
\sN \ot \sum_{1\dots4} \Oh(-d_i) \ot \sum_{j\ne 1} \Oh(-a_j) \ot
\Oh(-\si_1) \ot 0
\label{eq!J4}
\end{equation}
where $d_i=\wt x_i$, $a_j=\wt\Pf_j$ and $\si_1=\wt\Si_1$, and leave the
reader to think of names for the maps. Therefore $\sN$ has total Chern
class
\begin{equation}
\prod_{i=1}^4 (1-d_ih) \times (1-\si_1 h) \Big/ \prod_{j\ne1} (1-a_j h)
\label{eq!c2}
\end{equation}
The number of nodes $c_2(\sN)$ is then the $h^2$ term in the expansion
of \eqref{eq!c2}; recall that we view $h=c_1(\Oh_D(1))$ as an orbifold
class, so that $h^2=1/ab$ for $D=\PP(1,a,b)$.

\paragraph{Jer$_{12}$} The pivot $m_{12}$ appears in three Pfaffians
$\Pf_i=\Pf_{12,jk}$ for $\{i,j,k\}=\{3,4,5\}$ as the term
$m_{12}m_{jk}$, together with two other terms $m_{1j}m_{2k}$ of order
$\ge2$ in $x_{1\dots4}$. The Jacobian matrix restricted to $D$ thus has
three corresponding rows that are $m_{jk}$ times the same vector
$\partial m_{12}/\partial x_{1\dots4}$. This proportionality gives three
syzygies $\Si_l$ between these three rows, yoked by a second syzygy $T$
in degree $t=\hbox{adjunction number}-\wt m_{12}$. In other words, the
conormal bundle has the resolution
\begin{equation}
\sN\ot\bigoplus_4\Oh(-d_i)\ot\bigoplus_5\Oh(-a_j)\ot
\bigoplus_3\Oh(-\si_l) \ot\Oh(-t)\ot 0,
\end{equation}
so that the total Chern class of $\sN$ is the alternate product
\begin{equation}
\frac{\prod_4(1-d_ih) \prod_3(1-\si_lh)}{\prod_5(1-a_jh) (1-th)},
\end{equation}
with $c_2(\sN)$ equal to the $h^2$ term in this expansion.

\begin{exa} \label{rk!2v3} \rm
We read the number of nodes mechanically from the Hilbert numerator, the
matrix of weights and the choice of format. As a baby example, the
``interior'' projections of the two del Pezzo 3-folds of degree 6
discussed in \ref{s!same} have 2 and 3 respective nodes. These numbers
are the coefficient of $h^2$ in the formal power series
\begin{equation}
\frac{(1-h)^4(1-3h)}{(1-2h)^4}=1+h+2h^2 \enspace\hbox{and}\enspace
\frac{(1-h^4)(1-3h)^3}{(1-2h)^5(1-4h)}=1+h+3h^2.
\end{equation}

As a somewhat more strenuous example, in \eqref{eq!3.2},
\begin{description}

\item{Tom$_1$} has $\wt x_{1\dots4}=3,4,4,5$,
$\wt\Pf_{2\dots5}=8,8,7,6$, $\Si_1=10$, so that
\[
c(\sN) = \frac{\prod_{a\in [3,4,4,5,10]}(1-ah)}
{\prod_{b\in [6,7,8,8]}(1-bh)}=1+3h+28h^2,
\hbox{ giving $\frac{28}{1\cdot1\cdot2}=14$ nodes.}
\]

\item{Jer$_{25}$} has the same $x_i$, $\Pf_{1\dots5}=9,8,8,7,6$, $\Si_l=10,11,12$,
adjunction number = 19, $\wt m_{25}=5$, so
$c(\sN) = \frac{\prod_{a \in [3,4,4,5,10,11,12]}(1-ah)}
{\prod_{b \in [6,7,8,8,9,14]}(1-bh)}=1+3h+34h^2$, giving
$\frac{34}{1\cdot1\cdot2}=17$ nodes.

\end{description}

Try the other cases in \eqref{eq!3.2}--\eqref{eq!3.4} as homework.

\end{exa}

\section{Computer code and the GRDB database}
\label{s!grdb}

A Big Table with the detailed results of the calculations proving
Theorem~\ref{th!main} is online at the Graded Ring Database webpage
\begin{quote}
\verb!http://grdb.lboro.ac.uk! \quad + \hbox{Downloads.}
\end{quote}
This website makes available computer code implementing our calculations
systematically, together with the Big Table they generate. The code is
for the Magma system \cite{Ma}, and installation instructions are
provided; at heart, it only uses primary elements of any computer algebra
system, such as poly\-nomial ideal calculations and matrix
manipulations. The code runs online in the Magma Calculator
\begin{quote}
 \verb!http://magma.maths.usyd.edu.au/calc!
\end{quote}
All the data on the codimension~4 Fano 3-folds we construct is available
on webloc.\ cit.: follow the link to Fano 3-folds, select Fano index
$f=1$ (the default value), codimension~$=4$ and Yes for Projections of
Type~I, then submit. The result is data on the 116 Fano 3-folds with a
Type~I projection (the 116th is an initial case with $7\times12$
resolution, that projects to the complete intersection
$Y_{2,2,2}\subset\PP^6$ containing a plane, so is not part of our story
here). The $+$ link reveals additional data on each Fano.

The computer code follows closely the algorithm outlined as the proof of
Theorem~\ref{th!main}. For each Tom and Jerry format, we build a matrix
with random entries; some of these can be chosen to be single variables,
since we assume $Y$ is general for its format. We use row and column
operations to simplify the matrix further without changing the format.
The first failure tests (fishy zeroes, cone points and points of
embedding dimension~6) are now easy, and inspection of the equations on
affine patches at coordinate points on $Y$ is enough to determine
whether their local quotient weights are those of terminal
singularities. An ideal inclusion test checks that the singularities lie
on $D$. By good fortune, in every case that passes the tests so far, the
singular locus lies on one standard affine patch of $D$. We pass to this
affine patch and check that $\sI_{\Sing Y}\cdot\Oh_D$ defines a reduced
scheme there. We calculate the length of the quotient $\Oh_D/(\sI_{\Sing
Y}\cdot\Oh_D)$ on this patch, providing an alternative to the
computation of Section~\ref{s!Jac} (and a comforting sanity check).

The random entries in the matrix are not an issue: our non\-singularity
requirements are open, so if one choice leads to a successful $D\subset
Y$, any general choice also works. The only concern is false negative
reports, for example, an alleged nonreduced singular locus on $D$. To
tackle such hiccups, if a candidate fails at this stage (in practice, a
rare occurrence), we simply rerun the code with a new random matrix; the
fact that the code happens to terminate justifies the proof.

The conclusion is that every possible Tom and Jerry format for every
numerical Type~I projection either fails one of the human-readable tests
of Section~\ref{s!fail} (and we have made any number of such hand
calculations), or is shown to work by constructing a specific example.

To complete the proof of Theorem~\ref{th!main}, we check that the final
output satisfies the following two properties:
\begin{enumerate}
\renewcommand{\labelenumi}{(\alph{enumi})}

\item Every numerical candidate admits at least one Tom and one Jerry
unprojection.

\item Whenever a candidate has more than one Type~I centre, the
successful Tom and Jerry unprojections of any two correspond one-to-one,
with compatible numbers of nodes: the difference in Euler number
computed by the nodes is the same whichever centre we calculate from;
compare \eqref{eq!3.2}--\eqref{eq!3.4}.

\end{enumerate}

The polynomial ideal calculations of Nonsingularity
analysis~\ref{s!sings} (that is, the inclusion $I_D\subset\Rad(I_{\Sing
Y})$ and the statement that $\sI_{\Sing Y}\cdot\Oh_D$ is reduced) are
the only points where we use computer power seriously (other than to
handle hundreds of repetitive calculations accurately). In cases with 2
or 3 centres, even this could be eliminated by projecting to a complete
intersection and applying Bertini's theorem, as in
Section~\ref{exa!main}.

\section{Codimension~4 Gorenstein formats}\label{s!fmt}

The Segre embeddings $\PP^2\times\PP^2\subset\PP^8$ and
$\PP^1\times\PP^1\times\PP^1\subset\PP^7$ are well known codimension~4
projectively Gorenstein varieties with $9\times16$ resolution.
Singularity theorists consider the affine cones over them to be rigid,
because they have no nontrivial infinitesimal deformations or small
analytic deformation. Nevertheless, both are sections of higher
dimensional graded varieties in many different nontrivial ways. Each of
these constructions appears at many points in the study of algebraic
surfaces by graded rings methods.

\subsection{Parallel unprojection and extra\-symmetric format}
\label{s!xtra}

The extra\-symmetric $6\times6$ format occurs frequently, possibly first
in Dicks' thesis \cite{DPhD}. It is a particular case of triple
unprojection from a hypersurface in the product of three codimension 2
c.i.\ ideals. Start from the ``undeformed'' $6\times6$ skew matrix
\begin{equation}
M_0=
\begin{pmatrix}
b_3&-b_2&x_1&a_3&a_2\\
&b_1&a_3&x_2&a_1\\
&&a_2&a_1&x_3\\
&&&-b_3&b_2\\
&&&&-b_1
\end{pmatrix}
\end{equation}
with the ``extrasymmetric'' property that the top right $3\times3$ block
is symmetric, and the bottom right $3\times3$ block equals minus the top
left block. So instead of 15 independent entries, it has only 9
independent entries and 6 repeats.

Direct computation reveals that the $4\times4$ Pfaffians of $M_0$ fall
under the same numerics: of its 15 Pfaffians, 9 are independent and 6
repeats. One sees they generate the same ideal as the $2\times2$ minors
of the $3\times3$ matrix
\begin{equation}
N_0=
\begin{pmatrix}
x_1&a_3+b_3&a_2-b_2\\
a_3-b_3&x_2&a_1+b_1\\
a_2+b_2&a_1-b_1&x_3
\end{pmatrix}
\end{equation}
Here $N_0$ is the generic $3\times3$ matrix (written as symmetric plus skew),
with minors defining Segre $\PP^2\times\PP^2$, and thus far we have not
gained anything, beyond representing $\PP^2\times\PP^2$ as a nongeneric
section of $\Grass(2,6)$.

However $M_0$ can be modified to preserve the codimension~4 Gorenstein
property while destroying the sporadic coincidence with
$\PP^2\times\PP^2$. The primitive one-parameter way of doing this is to
choose the triangle $(1,2,6)$ and multiply the entries
$m_{12},m_{16},m_{26}$ by a constant $r_3$. This gives
\begin{equation}
M_1=
\begin{pmatrix}
r_3b_3&-b_2&x_1&a_3&r_3a_2\\
&b_1&a_3&x_2&r_3a_1\\
&&a_2&a_1&x_3\\
&&&-b_3&b_2\\
&&&&-b_1
\end{pmatrix}
\end{equation}
One checks that the three Pfaffians $\Pf_{12.i6}$ for $i=3,4,5$ are
$r_3$ times others, whereas three other repetitions remain unchanged. So
the $4\times4$ Pfaffians of $M_1$ still defines a Gorenstein
codimension~4 subvariety with $9\times16$ resolution. We can view it as
the Tom$_3$ unprojection of the codimension~3 Pfaffian ideal obtained by
deleting the final column, with $x_3$ as unprojection variable.

If $r_3=\rho^2$ is a perfect square then floating the square root $\rho$ to
the complementary entries $m_{34},m_{35},m_{45}$ restores the original
extrasymmetry. In general this is a ``twisted form'' of $\PP^2\times\PP^2$:
changing the sign of $\rho$ swaps the two factors.

A more elaborate version of this depends on 8 parameters:
\begin{equation}
M_2=
\begin{pmatrix}
r_3s_0b_3&-r_2s_0b_2&x_1&r_2s_1a_3&r_3s_1a_2&\\
&r_1s_0b_1&r_1s_2a_3&x_2&r_3s_2a_1&\\
&&r_1s_3a_2&r_2s_3a_1&x_3&\\
&&&-r_0s_3b_3&r_0s_2b_2&\\
&&&&-r_0s_1b_1&
\end{pmatrix}
\label{eq!r1r2r3}
\end{equation}

Now the same three Pfaffians $\Pf_{12.i6}$ are divisible by
$r_3$, and the complementary three are divisible by $s_3$ with
the same quotient, so one has to do a little cancellation to
see the irreducible component. The necessity of cancelling
these terms (although cheap in computer algebra as the colon
ideal) has been a headache in the theory for decades, since it
introduces apparent uncertainty as to the generators of the ideal.

The right way to view this is as the triple parallel unprojection of
the hypersurface
\begin{equation}
V(a_1a_2b_3r_3s_3 + a_1a_3b_2r_2s_2 + a_2a_3b_1r_1s_1 + b_1b_2b_3r_0s_0)
\end{equation}
in the product ideal $\prod_{i=1}^3(a_i,b_i)$. Then
\[
x_1=\frac{a_2a_3r_1s_1+b_2b_3r_0s_0}{a_1}=
-\frac{a_2b_3r_2s_2+a_2b_3r_3s_3}{b_1},
\]
etc., and the ideal is generated by the Pfaffians of the three matrices
\begin{equation}
\left(
\begin{smallmatrix}
x_2&b_1r_0s_0&a_1r_3s_3&a_3 \\
&-a_1r_2s_2&-b_1r_1s_1&b_3 \\
&&x_3&a_2 \\
&&&b_2
\end{smallmatrix}\right), \
\left(
\begin{smallmatrix}
x_1&b_3r_0s_0&a_3r_2s_2&a_2 \\
&-a_3r_1s_1&-b_3r_3s_3&b_2 \\
&&x_2&a_1 \\
&&&b_1
\end{smallmatrix}\right), \
\left(
\begin{smallmatrix}
x_3&b_2r_0s_0&a_2r_1s_1&a_1 \\
&-a_2r_3s_3&-b_2r_2s_2&b_1 \\
&&x_1&a_3 \\
&&&b_3
\end{smallmatrix}\right).
\notag
\end{equation}
If the $r_i$ and $s_i$ are nonzero constants, one still needs the
square root of the discriminant $\prod_{i=0}^3 (r_is_i)$ to get back
to $\PP^2\times\PP^2$.

\subsection{Double Jerry}\label{s!nJ}
The equations of Segre $\PP^1\times\PP^1\times\PP^1\subset\PP^7$ are the
minors of a $2\times2\times2$ array; they admit several extensions, and
it seems most likely that there is no irreducible family containing them
all. One family consists of various ``rolling factors'' formats
discussed below; here we treat ``double Jerry''.

Start from the equations written as
\begin{equation}
\begin{array}{cl}
sy_i = x_jx_k & \hbox{for $\{i,j,k\}=\{1,2,3\}$}, \\
tx_i = y_jy_k & \hbox{for $\{i,j,k\}=\{1,2,3\}$}, \\
st = x_iy_i & \hbox{for $i=1,2,3$}.
\end{array}
\label{eq!cube}
\end{equation}
corresponding to a hexagonal view of the cube centred at vertex $s$
(with three square faces $\square sx_iy_kx_j$, and $t$ behind the
page, cf.~\eqref{eq!cube1}):
\begin{equation}
\renewcommand{\arraycolsep}{.125em}
\begin{matrix}
&&y_2 \\[-3pt]
x_3&&&&&x_1 \\[-3pt]
&& s \\[-5pt]
y_1&&&&&y_3 \\[-3pt]
&&x_2
\end{matrix}
\begin{picture}(0,0)(25,0)
\qbezier(-5,19)(-8,17)(-11,15)
\qbezier(5,19)(8,17)(11,15)
\qbezier(-5,-15)(-8,-13)(-11,-11)
\qbezier(5,-15)(8,-13)(11,-11)
\qbezier(-5,3)(-8,5)(-11,7)
\qbezier(5,3)(8,5)(11,7)
\qbezier(0,-4)(0,-8)(0,-12)
\qbezier(-17,6)(-17,2)(-17,-2)
\qbezier(0,-5)(0,-9)(0,-13)
\qbezier(17,6)(17,2)(17,-2)
\end{picture}
\end{equation}
Eliminating both $s$ and $t$ gives the codimension~2 c.i.
\begin{equation}
(x_1y_1 = x_2y_2 = x_3y_3) \subset\PP^5,
\label{eq!xiyi}
\end{equation}
containing the two codimension~3 c.i.s $\bx=0$ and $\by=0$ as divisors.
We can view $\bx$ as a row vector and $\by$ a column vector, and the
two equations \eqref{eq!xiyi} as the matrix products
\begin{equation}
\bx A\by = \bx B\by = 0, \quad\hbox{where}\quad
A=\left(\begin{smallmatrix}
1&0&0 \\
0&-1&0 \\
0&0&0
\end{smallmatrix}\right)\!,\,
B=\left(\begin{smallmatrix}
0&0&0 \\
0&1&0 \\
0&0&-1
\end{smallmatrix}\right)\!.
\end{equation}
The unprojection equations for $s$ and $t$ separately take the form
\begin{equation}
t\bx = (A\by) \times (B\by)  \quad\hbox{and}\quad
s\by = (\bx A) \times (\bx B),
\end{equation}
where $\times$ is cross product of vectors in $\C^3$, with the
convention that the cross product of two row vectors is a column vector
and vice versa. For example, $\bx A=(x_1,-x_2,0)$, $\bx B=(0,x_2,-x_3)$
and the equations $s\by=(\bx A) \times (\bx B)$ giving the first line of
\eqref{eq!cube} are deduced via Cramer's rule from \eqref{eq!xiyi}.

We can generalise this at a stroke to $A,B$ general $3\times3$ matrices.
That is, for $\bx$ a row vector and $\by$ a column vector, $\bx A\by=\bx
B\by=0$ is a codimension~2 c.i.; since these are general bilinear forms
in $\bx$ and $\by$, it represents a universal solution to two elements
of the product ideal $(x_1,x_2,x_3)\cdot(y_1,y_2,y_3)$. It has two
single unprojections:
\begin{align}
\bx A\by=\bx B\by=0 \quad\hbox{and}\quad
s\by=(\bx A) \times (\bx B),
\label{eq!sby} \\
\bx A\by=\bx B\by=0 \quad\hbox{and}\quad
t\bx=(A\by) \times (B\by),
\label{eq!tbx}
\end{align}
either of which is a conventional $5\times5$ Pfaffian, and a parallel
unprojection putting those equations together with a 9th {\em long
equation}
\begin{equation}
st = \hbox{something complicated}.
\end{equation}
The equation certainly exists by the Kustin--Miller theorem. It can be
obtained easily in computer algebra by coloning out any of
$x_1,x_2,x_3,y_1,y_2,y_3$ from the ideal generated by the eight
equations \eqref{eq!sby} and \eqref{eq!tbx}. Its somewhat amazing right
hand side has 144 terms, each bilinear in $x,y$ and biquadratic in
$A,B$. Taking a hint from $144=12\times12$, we suspect that it may have
a product structure of the form
\begin{equation}
\bx(A \wedge B)\times(A\wedge B)\by,
\end{equation}
with ``$\times$'' and ``$\wedge$'' still requiring elucidation.

If the entries of $A$ and $B$ are constants, one gets back to
$\PP^1\times\PP^1\times\PP^1$ after coordinate changes based on
the three roots $(\la_i:\mu_i)$ of the relative characteristic
equation $\det(\la A-\mu B)=0$ and the three eigenvectors
$v_i=\ker(\la_i A-\mu_i B)$. Swapping the roots permutes the
three factors.

The significance of the double Jerry parallel unprojection format is
that it covers any Jerry case where the pivot is one of the generators
of $I_D$. Indeed, if the regular sequence generating $I_D$ is
$s,x_1,x_2,x_3$, a Jerry matrix for $D$ is
\begin{equation}
\begin{pmatrix}
s&m_{13}&m_{14}&m_{15} \\
&m_{23}&m_{24}&m_{25} \\
&& y_3 & -y_2 \\
&&& y_1
\end{pmatrix} \quad\hbox{where}\quad
\begin{array}{l}
(m_{13},m_{14},m_{15})=\bx A, \\[3pt]
(m_{23},m_{14},m_{15})=\bx B.
\end{array}
\end{equation}
for some $3\times3$ matrices $A,B$. Unprojecting $D$ gives a double Jerry.

\subsection{Rolling factors format}\label{s!rf}

Rolling factors view a divisor $X\subset V$ on a normal projective
variety $V\subset\PP^n$ as residual to a nice linear system. This
phenomenon occurs throughout the literature, with typical cases a
divisor on the Segre embedding of $\PP^1\times\PP^3$, or on a rational
normal scroll $\FF$, or on a cone over a Veronese embedding. A divisor
$X\subset\PP^1\times\PP^3$ in the linear system
$|ah_1+(a+2)h_2|=|{-}K_V+bH|$ is of course defined by a single
bihomogeneous equation in the Cox ring of $\PP^1\times\PP^3$, but to get
equations in the homogeneous coordinate ring of Segre
$\PP^1\times\PP^3\subset\PP^7$ we have to add $|2h_1|$. This is a type
of hyperquotient, given by one equation in a nontrivial eigenspace.

Dicks' thesis \cite{DPhD} discussed the generic pseudoformat
\begin{equation}
\begin{gathered}
\bigwedge^2 \begin{pmatrix} a_1&a_2&a_3&a_4 \\ b_1&b_2&b_3&b_4
\end{pmatrix}=0, \quad\hbox{and} \\[4pt]
\renewcommand{\arraycolsep}{.25em}
\begin{array}{rclcc}
m_1a_1+m_2a_2+m_3a_3+m_4a_4 &=& 0 \\
m_1b_1+m_2b_2+m_3b_3+m_4b_4 &\equiv&
n_1a_1+n_2a_2+n_3a_3+n_4a_4 &=& 0 \\
&& n_1b_1+n_2b_2+n_3b_3+n_4b_4 &=& 0.
\end{array}
\label{eq!roll}
\end{gathered}
\end{equation}
One sees that under fairly general assumptions the ``scroll'' $V$
defined by the first set of equations of \eqref{eq!roll} is
codimension~3 and Cohen--Macaulay, with resolution
\[
\Oh_V\ot R\ot 6R\ot 8R\ot 3R\ot0.
\]

On the right, the identity is a preliminary condition on quantities in
the ambient ring. If we assume (say) that $R$ is a regular local ring
and $a_i,b_i,m_i,n_i\in R$ satisfy it (and are ``fairly general''), the
second set defines an elephant $X\in|{-}K_V|$ (anticanonical divisor)
which is a codimension~4 Gorenstein variety with $9\times16$ resolution.

The identity in \eqref{eq!roll} is a quadric of rank~16. It is a little
close-up view of the ``variety of complexes'' discussed in \cite{Ki},
Section~10. To use this method to build genuine examples, we have to
decide how to map a regular ambient scheme into this quadric; there are
several different solutions. If we take the $a_i,b_i$ to be independent
indeterminates, the first set of equations gives the cone on Segre
$\PP^1\times\PP^3\subset\PP^7$, and the second set consists of a single
quadratic form $q$ in 4 variables evaluated on the two rows, so that
$X\subset V$ is given by $q(\ba)=\fie(\ba,\bb)=q(\bb)=0$, with $\fie$
the associated symmetric bilinear form (cf.\ \eqref{eq!rf}). This format
seems to be the only commonly occurring codimension~4 Gorenstein format
that tends not to have any Type~I projection.

On the other hand, if there are coincidences between the $a_i,b_i$,
there may be other ways of choosing the $m_i,n_i$ to satisfy the
identity in \eqref{eq!roll} without the need to take $m_i,n_i$ quadratic
in the $a_i,b_i$: for example, if $a_2=b_1$, we can roll $a_1\to a_2$
and $b_1\to b_2$.

\bigskip
\noindent
Gavin Brown \\
School of Mathematics,
Loughborough University \\
LE11 3TU, UK \\
\verb!G.D.Brown@lboro.ac.uk!

\medskip
\noindent
Michael Kerber\\
Institute of Science and Technology (IST) Austria\\
3400 Klosterneuburg, Austria\\
\verb!michael.kerber@ist.ac.at!

\medskip
\noindent
Miles Reid\\
Mathematics Institute,
University of Warwick\\
Coventry CV4 7AL, UK\\
\verb!Miles.Reid@warwick.ac.uk!

\end{document}